\documentclass[journal, subeqn, caption]{IEEEtran}
\usepackage[singlelinecheck=false,labelsep=period,font=footnotesize]{caption}
\usepackage{comment}
\usepackage{mathtools}
\usepackage{cite}
\usepackage{float}
\usepackage{multicol}
\usepackage{multirow}
\usepackage{hyperref}
\usepackage{graphicx}
\usepackage{amsfonts}
\usepackage{amsmath}
\usepackage{epstopdf}
\usepackage{subcaption}
\usepackage{eurosym}
\usepackage{enumitem}
\usepackage{tikz,pgfplots,pgfplotstable,booktabs}
\pgfplotsset{compat=1.8}
\usepgfplotslibrary{statistics}
\usepackage{balance}
\newenvironment{ldescription}[1]
  {\begin{list}{}%
   {\renewcommand\makelabel[1]{##1\hfill}%
   \settowidth\labelwidth{\makelabel{#1}}%
   \setlength\leftmargin{\labelwidth}
   \addtolength\leftmargin{\labelsep}}}
  {\end{list}}

\begin{document}
\title{Feature-driven Improvement of Renewable Energy Forecasting and Trading}

\author{M. A. Mu\~noz, J. M. Morales, S. Pineda

\thanks{M. A. Mu\~noz, J. M. Morales (corresponding author) and S. Pineda are with the research group OASYS at the University of Malaga, Malaga, Spain. E-mails: miguelangeljmd@uma.es; juan.morales@uma.es; spinedamorente@gmail.com.% <-this % stops a space

This work was supported, in part, by the European
Research Council (ERC) under the EU Horizon 2020 research and innovation
programme (grant agreement No. 755705) and, in part, by the Spanish Ministry of Economy, Industry, and Competitiveness through project ENE2017-83775-P.}}

\vspace{-15mm}

\maketitle

\begin{abstract}
Inspired from recent insights into the common ground of machine learning, optimization and decision-making, this paper proposes an easy-to-implement, but effective procedure to enhance both the quality of renewable energy forecasts and the competitive edge of renewable energy producers in electricity markets with a dual-price settlement of imbalances. The quality and economic gains brought by the proposed procedure essentially stem from the utilization of valuable predictors (also known as \emph{features}) in a data-driven newsvendor model that renders a computationally inexpensive linear program. We illustrate the proposed procedure and numerically assess its benefits on a realistic case study that considers the aggregate wind power production in the Danish DK1 bidding zone as the variable to be predicted and traded. Within this context, our procedure leverages, among others, spatial information in the form of wind power forecasts issued by transmission system operators (TSO) in surrounding bidding zones and publicly available in online platforms. We show that our method is able to improve the quality of the wind power forecast issued by the Danish TSO by several percentage points (when measured in terms of the mean absolute or the root mean square error) and to significantly reduce the balancing costs incurred by the wind power producer.
%Wind power is not easily predictable which force wind farm owners to deal with uncertainties when bidding in electricity markets and iccur in costs do to imbalances. In this paper we show how wind power producers can benefit from using data driven strategies and exogenous information known as features to improve both their predictions and bidding strategies. A case study is presented in which DK1 Denmark bidding zone is considered as an aggregated wind farm. Wind power producer cost function is model as a newsvendor cost function extended to take into account additional information. Among others, we use as feature the forecast of surrounding bidding zones as spatial information. We propose a two-step method to obtain an offer close enough to the actual value but modified to take into account the historical cost information.
%
\end{abstract}

\vspace{-0.2cm}

\section*{Nomenclature}
\subsection{Sets and Indices}
\begin{ldescription}{$xxxx$}
\item [$j$] Index of features.
\item [$t$] Index of time periods.
\item[$\mathcal{T}$] Training set.
\item[$\widetilde{\mathcal{T}}$] Test set.
\end{ldescription}

\subsection{Parameters}
\begin{ldescription}{$xxxx$}
\item [$\overline{E}$] Maximum hourly wind energy production (MWh).
\item [$E_t$] Actual wind energy produced at hour $t$ (MWh).
\item [$\lambda^{B}_t$] Balancing market price at hour $t$ (\euro/MWh).
\item [$\lambda^{D}_t$] Day-ahead market price at hour $t$ (\euro/MWh).
\item [$\lambda^{-}_t$] Upward regulation price in the balancing market at hour~$t$ (\euro/MWh).
\item [$\lambda^{+}_t$] Downward regulation price in the balancing market at hour $t$ (\euro/MWh).
\item [$\psi^{-}_t$] Marginal opportunity cost for underproduction at hour $t$ (\euro/MWh).
\item [$\psi^{+}_t$] Marginal opportunity cost for overproduction at hour $t$ (\euro/MWh).
\end{ldescription}

\subsection{Variables}
\begin{ldescription}{$xxxx$}
\item [$E^D_t$] Energy bid for hour $t$ of the market horizon submitted to the day-ahead electricity market (MWh).
\end{ldescription}

\section{Introduction}
\label{sec:introduction}
Thrilling yet challenging times lie ahead for the electrical power industry. The development of microgrids, the growing contribution of weather-driven renewable energy sources, the higher involvement of power consumers, and the increasing exchange of electricity among neighbouring regions are demanding profound changes in the power sector. These changes are expected to turn power systems into complex and critical cyber-physical systems, where data will be generated and made accessible in abundance and where data will play an increasingly important role for decision-making.

In Europe, for example, the efforts invested by the EU member countries in setting up a single electricity market have been accompanied with the development of the so-called \emph{ENTSO-e Transparency Platform} \cite{entsoe}, a web database where data on electricity generation, transmission and consumption in the pan-European market is gradually collected, published, and made publicly available for download. In fact, the research here described constitutes an example of how the information gathered in this platform can be used to generate extra value in two important tasks that are performed daily in electricity markets, namely, renewable energy forecasting and trading. More specifically, we focus on the aggregate onshore wind power production of the DK1 bidding zone of the European market and show that the forecast that is issued by the Danish TSO everyday can be noticeably and easily improved by leveraging the information contained in that platform, in particular, the forecasts of wind power production in neighbouring areas issued by their respective TSOs. Furthermore, we also show that this very same information can be used to increase the profitability of wind power production in electricity markets with a dual-price financial settlement for imbalances.

To achieve these goals, we exploit recent insights into the close bonds that connect the fields of machine learning, optimization and decision-making. For some years now, researchers from these fields have been developing methods that leverage data not  to make better predictions, but to make better decisions, on the grounds that the former does not always necessarily imply the latter. In this line, we mention the works \cite{Bertsimas2014, VanParys2017a, Elmachtoub2017, Ban2018}. From among them, our work builds on the data-driven model for the newsvendor problem developed in \cite{Ban2018}, because of its simplicity and because it neatly fits with the setup of our problem. As explained later, however, our problem exhibits some peculiarities that make it especially challenging.

On a different front, the technical literature on wind power forecasting and trading is tremendously vast. Mentioning all the many relevant references on both topics in this paper would be, therefore, an infeasible and purposeless task. We refer, instead, to monographs \cite{kariniotakisBOOK, Morales2014book}, which offer a comprehensive treatment of both topics, and highlight next those works that, we believe, are most closely related to ours. In the realm of wind power prediction, such works would be those that either seek to model the spatial correlations among wind sites (see, e.g., \cite{MORALES2010, tastu2011, Xie2014, Zhao2018}) or to adaptively combine alternative wind power forecasts  for the same site so as to produce a better one (see, for example, \cite{Nielsen2007, SANCHEZ2006, SANCHEZ2008}). In our case, however, we do not aim at developing a better forecasting model. What we propose, instead, is a general mathematical framework to improve the forecasts delivered by any \emph{existing} method by leveraging available power system data. To do so, we use a straightforward procedure that exploits extra information, for example, information on spatially correlated phenomena. On the other hand, there also exists a wealth of methods to determine the optimal energy bid that a wind power producer should place in a day-ahead electricity market (see, for instance, \cite{Pinson2007, Morales2010b, Bitar2011, Zugno2013, Matevosyan2006, Bathurst2002, Dent2011}). To this end, all these methods make explicit use of stochastic models for the wind power production and/or market prices, for example, in the form of scenario forecasts or predictive densities.
Additionally, other strategies have been also proposed to cope with the inherent uncertainty in wind power production such as the purchase of power reserves \cite{Du2017} or by means of a combined portfolio of wind and hydro power generation \cite{Nieta2013}. What distinguishes our work from these others is that we directly derive a wind power day-ahead bid from available \emph{point} forecasts and other relevant data, thus avoiding the need to generate \emph{scenarios} or \emph{probabilistic} forecasts for electricity prices and wind power production.

%the ability of our approach to easily derive a wind power day-ahead bid from any available information on any variable potentially relevant to reduce the expected imbalance cost incurred by the wind power producer.

%We particularly mention the recently published paper \cite{Carriere2019} as the work that is probably closest to ours.

We particularly mention \cite{Mazzi2016} and the recently published paper \cite{Carriere2019} as the works that are probably closest to ours. In \cite{Mazzi2016}, a reinforcement learning algorithm is built to compute and follow the nominal level of the profit-maximizing quantile forecast of wind power that should be bid into the day-ahead market. While their algorithm is designed to learn and track the expected marginal opportunity costs directly from market data, they assume that a good estimate of the wind power predictive density is available (as in the other references mentioned above). Our approach, on the contrary, is freed from this classical assumption. In \cite{Carriere2019}, the authors propose two data-driven approaches to reduce the imbalance costs incurred by renewable energy producers. In their first approach, they formulate a meta-optimization problem whereby the hyper-parameters of all the forecasting models involved in the decision-making process are tuned to minimize the imbalance costs. In their second approach, they directly train an artificial neural network to that very same end. In contrast with our proposal, which boils down to a linear programming problem, the complexity of theirs is such that they need to resort to heuristic optimization algorithms. Furthermore, our way to produce market bids is somewhat different: We do not search for a bidding model that overrides the need for forecasts (understood in a classical statistical sense), but collect those forecasts, among other features, and combine them to produce profit-maximizing bids.
%Furthermore, our proposal is not to train the forecasting models to minimize the imbalance costs or directly to get rid of those models, but to use all the information available (also in the form of forecasts) to produce a more cost-effective renewable power bid.

% The contributions of our paper are, therefore, the following.
% %
% \begin{enumerate}
%     \item We propose a method to improve a forecast of renewable power production by leveraging extra information on potentially correlated phenomena, such as the forecasts of the renewable power production in adjacent regions. The method is based on a data-driven model for the newsvendor problem.
%     \item We introduce a variant of the method proposed in point 1) above to increase the profitability of renewable power production in electricity markets with a dual-price settlement for imbalances.
%     \item We illustrate the benefits of our approach on a realistic case study that considers the aggregate onshore wind power production of the DK1 bidding area of the European market.
% \end{enumerate}
%

The contributions of our paper are, therefore, the following.
\begin{enumerate}
    \item We propose a general data-driven optimization framework to improve a forecast of renewable power production by
    \begin{enumerate}
        \item tailoring the optimization problem by which the forecast is improved to the specific use for which the forecast is intended, and
        \item leveraging extra information on potentially related phenomena, such as the forecasts of the renewable power production in adjacent regions.
    \end{enumerate}
    The proposed framework is based on a data-driven model for the newsvendor problem and reduces to a computationally efficient and easy-to-implement linear program.
    \item We focus on two particular instances of the general framework proposed in point 1) above. The first one seeks to improve the estimate of the median of the renewable energy production, while the second aims to increase the profitability of renewable power production in electricity markets with a dual-price settlement for imbalances.
    \item We illustrate the benefits of our approach on a realistic case study that considers the aggregate onshore wind power production of the DK1 bidding area of the European market. Furthermore, this case study serves us to argue that the accuracy of the renewable energy forecasts issued by TSOs could be significantly improved, if they were willing to share those forecasts among them.
\end{enumerate}

The rest of this paper is organized as follows. Section \ref{sec:methodology} states the problem and introduces the data-driven newsvendor model we propose to solve it, while Section \ref{sec:design&training} elaborates on its practical implementation. Numerical results from the application of our approach to real data are discussed in Section~\ref{sec:results}. Finally, conclusions are duly drawn in Section~\ref{sec:conclusions}.

\section{Data-driven Approach}
\label{sec:methodology}
Consider an electricity market for short-term energy transactions that consists of a day-ahead market and a \emph{dual-price} balancing market. In the former, energy offers and bids are typically to be submitted between 12 and 36 hours in advance of the actual delivery of electricity. In the latter, deviations of market participants with respect to their day-ahead dispatch are financially settled at a price that depends on the sign of the total system imbalance \cite[Ch. 7]{Morales2014book}.

%In such a context, the optimal offer $E^{D}$ that a (price-taker) risk-neutral renewable energy producer should place in the day-ahead market is given as the solution to the following linear programming problem, whereby the renewable energy producer seeks to minimize the expected opportunity costs for under- and overproduction:
%
% \begin{align}\label{eq:nv}
%     \min_{E^D \in [0, \overline{E}]} \enskip  \mathbb{E}\left[\psi^{-}(E^D - E)^{+} + \psi^{+}(E - E^D)^{+}\right]
% \end{align}
% %
% where $(x)^{+}:= \max(x,0)$.

% In problem~\eqref{eq:nv}, the expectation is taken over the stochastic input parameters $E$, $\psi^{-}$ and $\psi^{+}$. These parameters represent the renewable energy production and the marginal opportunity costs for under- and overproduction, respectively. Logically, these parameters are uncertain to the renewable energy producer at the moment of offering in the day-ahead market, and as such, the way the solution to problem~\eqref{eq:nv} is addressed depends on the information on $E$, $\psi^{-}$ and $\psi^{+}$ we have.

In such a context, the market revenue $\rho$ of a renewable energy producer in a dual-price balancing market is given by
\begin{equation}\label{eq:revenue}
        \rho = \lambda^D E - \left( \psi^{-}(E^D - E)^{+} + \psi^{+}(E - E^D)^{+} \right)
\end{equation}
where $\lambda^D$, $E^{D}$, $\psi^{-}$, $\psi^{+}$, and $E$ represent the day-ahead market price, the day-ahead renewable energy bid, the marginal opportunity costs for under- and overproduction, and the eventual renewable energy production, respectively. In~\eqref{eq:revenue}, the first term accounts for the incomes the renewable power producer would obtain from partaking in the day-ahead market if she had perfect information on her eventual production, while the second is the opportunity cost the producer incurs in deviating from the day-ahead bid $E^{D}$.
Logically, parameters  $\psi^{-}$, $\psi^{+}$, and $E$ are uncertain to the renewable energy producer at the moment of offering in the day-ahead market. Besides, the term $\lambda^D E$ is out of the power producer's control. As a result, the optimal offer $E^{D}$ that a (price-taker) risk-neutral renewable energy producer should place in the day-ahead market is given as the solution to the following linear programming problem, whereby the renewable energy producer seeks to minimize the expected opportunity cost for under- and overproduction:
\begin{align}\label{eq:nv}
    \min_{E^D \in [0, \overline{E}]} \enskip  \mathbb{E}\left[\psi^{-}(E^D - E)^{+} + \psi^{+}(E - E^D)^{+}\right]
\end{align}
where $(x)^{+}:= \max(x,0)$.

In problem~\eqref{eq:nv}, the expectation is taken over the stochastic input parameters $E$, $\psi^{-}$ and $\psi^{+}$. Actually, the way the solution to problem~\eqref{eq:nv} is addressed depends on the information we have about these parameters. Furthermore, this problem must be (independently) solved for every trading period comprising the day-ahead market horizon (typically the 24 hours of a day). For simplicity, though, we have dropped the time index from the problem formulation. We will introduce that index in a later stage of our exposition.

The marginal opportunity costs for under- and overproduction, i.e., $\psi^{-}$ and $\psi^{+}$, are defined as:
\begin{align}
    \psi^{-} = \lambda^{-} - \lambda^D \label{eq:psi_def}\\
    \psi^{+} = \lambda^D  - \lambda^{+}
\end{align}
where, in turn, the prices for under- and overproduction, i.e., $\lambda^{-}$ and $\lambda^{+}$ are given by:
\begin{equation}  \label{eq:priceup}
    \lambda^{-} =
        \begin{cases}
          \lambda^{B} & \text{if } \lambda^{B} \ge \lambda^{D} \\
          \lambda^{D} & \text{if } \lambda^{B} < \lambda^{D} \\
        \end{cases}
\end{equation}
\begin{equation}  \label{eq:pricedw}
    \lambda^{+} =
        \begin{cases}
          \lambda^{D} & \text{if } \lambda^{B} \ge \lambda^{D} \\
          \lambda^{B} & \text{if } \lambda^{B} < \lambda^{D} \\
       \end{cases}
\end{equation}
In \eqref{eq:priceup} and \eqref{eq:pricedw}, $\lambda^{D}$ and $\lambda^{B}$ denote the day-ahead and the balancing market prices, in that order.

Therefore, according to the rules \eqref{eq:psi_def}--\eqref{eq:pricedw} of a dual-price imbalance settlement, the overproduction of a renewable energy producer is always rewarded at a price lower than or equal to the day-ahead market price, while their underproduction is always penalized at a price higher than or equal to the day-ahead market price. This settlement is, at least, used in some European countries such as Spain and Denmark \cite{wang2015}.

Problem~\eqref{eq:nv} takes the form of the classical \emph{newsvendor} problem
\cite{qin2011}, for which an analytical solution exists. Indeed, the optimal solution to this problem (that is, the optimal bid $E^{D*}$), is given by:
\begin{equation}\label{eq:ana_sol}
    E^{D*} = F^{-1}_{E}\left(\frac{\bar{\psi}^{+}}{\bar{\psi}^{+}+\bar{\psi}^{-}}\right)
\end{equation}
where $F_{E}$ is the cumulative distribution function (cdf) of the renewable energy production corresponding to the time period of the market horizon for which the day-ahead bid must be submitted, and the overbar character denotes the \emph{expected value} of the random variable underneath.

Despite its apparent simplicity, the application of formula \eqref{eq:ana_sol} is quite demanding, as it requires models to produce a probabilistic forecast of $E$ (i.e., an estimate of its cdf) and point forecasts of $\psi^{-}$ and $\psi^{+}$. In the first approach proposed in \cite{Carriere2019}, for example, those models are tuned (by way of what they call a meta-optimization problem) to produce a good estimate of \eqref{eq:ana_sol}. Our goal, though, is to sidestep the need for those models and directly use available data instead. This motivates our data-driven approach, which we gradually build next.

Suppose that the renewable energy producer is to place a bid in the day-ahead market and that measurements of her renewable energy production at past periods are available. We can then directly use the \emph{empirical} cdf of these data, namely, $\widehat{F}_{E}$, in lieu of $F_{E}$ in \eqref{eq:ana_sol}, which thus becomes
\begin{align} \label{eq:q_opt_emp}
    {\widehat{E}}^{D} = \inf \left\{y: \widehat{F}_{E}(y) \ge \frac{\bar{\psi}^{+}} {\bar{\psi}^{+} + \bar{\psi}^{-}}\right\}
\end{align}
where the infimum is required due to the discrete nature of $\widehat{F}_{E}$. Naturally, ${\widehat{E}}^{D}$ in \eqref{eq:q_opt_emp} and $E^{D*}$ in \eqref{eq:ana_sol} are generally different, and therefore, ${\widehat{E}}^{D}$ is usually suboptimal in~\eqref{eq:nv}. Actually, ${\widehat{E}}^{D}$ is the  solution to the following \emph{sample average approximation} (SAA) of \eqref{eq:nv}
\begin{align}\label{eq:nv_dis}
    \min_{E^D \in [0, \overline{E}]}  \frac{1}{|\mathcal{T}|} \sum_{t\in\mathcal{T}} \bar{\psi}^{-}(E^D - E_t)^{+} + \bar{\psi}^{+}(E_t - E^D)^{+}
\end{align}
From the equivalence between~\eqref{eq:q_opt_emp} and \eqref{eq:nv_dis}, we can infer that if we (artificially) set $\bar{\psi}^{-}=\bar{\psi}^{+}=1$ in \eqref{eq:nv_dis}, we get an estimate of the \emph{median} of the renewable energy production. We will leverage this fact later on to develop a straightforward method to enhance the quality of renewable energy forecasts.

Problem~\eqref{eq:nv_dis}, however, is likely to deliver poor bids ${\widehat{E}}^{D}$, because it overlooks the fact that, at the moment of bidding, the renewable power producer may have information on a vector $\mathbf{x}$ of $p$ features with some predictive power on her future production. Accordingly, to get a better bid ${\widehat{E}}^{D}$, we need to reformulate the SAA problem~\eqref{eq:nv_dis} to account for and take advantage of that information. For this purpose, we consider the enriched dataset $\left\{(E_t, \mathbf{x}_t),\forall t \in \mathcal{T}\right\}$, where $\mathbf{x}_t$ is the $p$-dimensional realization of features $\mathbf{x}$ observed at time $t$. These features may include measures of potentially explanatory variables available at time period $t$ or forecasts of these variables issued for that time period. We then follow the approach proposed in \cite{Ban2018} and consider the following linear decision rule
\begin{align} \label{eq:linear-rule}
    \mathcal{Q}=\Big \{ E^D: \mathcal{X} \to \mathbb{R}:E^D(x) = \mathbf{q \cdot x} = \sum_{j=1}^p q^j x^j \Big \},
\end{align}
which, inserted into~\eqref{eq:nv_dis}, renders
\begin{align}
    &\min_{\mathbf{q}}  \frac{1}{|\mathcal{T}|} \sum_{t\in\mathcal{T}} \bar{\psi}^{-}\left(\sum_{j=1}^p q^j x^j_t - E_t\right)^{+} + \bar{\psi}^{+}\left(E_t - \sum_{j=1}^p q^j x^j_t\right)^{+}\label{eq:nv_dis_feature1}\\
    & \text{s. t.} \enskip 0 \leq \sum_{j=1}^p q^j x^j_t \leq \overline{E}, \enskip \forall t\in\mathcal{T} \label{eq:nv_dis_feature2}
\end{align}
Nonetheless, problem~\eqref{eq:nv_dis_feature1}--\eqref{eq:nv_dis_feature2} still requires further elaboration to become a fully data-driven model. Indeed, while in the technical literature on the data-driven newsvendor problem (see, for instance, \cite{Ban2018} and \cite{HUBER2019}), the marginal opportunity costs $\bar{\psi}^{-}$ and $\bar{\psi}^{+}$ are assumed to be known with certainty, in our case, these costs are unknown to the renewable energy producer at the moment of bidding into the day-ahead market. Consequently, problem~\eqref{eq:nv_dis_feature1}--\eqref{eq:nv_dis_feature2} still needs the support of a forecasting model that provides it with an estimate of $\bar{\psi}^{-}$ and $\bar{\psi}^{+}$. To circumvent this hurdle, we propose to work with the even more enriched dataset $\left\{(E_t, \psi_t^{-}, \psi_t^{+},  \mathbf{x}_t), \forall t\in\mathcal{T} \right\}$, where the pair $(\psi_t^{-}, \psi_t^{+})$ represents the marginal costs of under- and overproduction that were observed at time $t$, and solve instead the following optimization problem:
\begin{align}
    &\min_{\mathbf{q}}  \frac{1}{|\mathcal{T}|} \sum_{t\in\mathcal{T}} \psi_t^{-}\left(\sum_{j=1}^p q^j x^j_t - E_t\right)^{+}\! +\! \psi_t^{+}\left(E_t - \sum_{j=1}^p q^j x^j_t\right)^{+}\label{eq:nv_dis_feature_dd1}\\
    & \text{s. t.} \enskip 0 \leq \sum_{j=1}^p q^j x^j_t \leq \overline{E}, \enskip \forall t\in\mathcal{T} \label{eq:nv_dis_feature_dd2}
\end{align}
where we have replaced $\bar{\psi}^{-}$ and $\bar{\psi}^{+}$ with $\psi_t^{-}$ and $\psi_t^{+}$, respectively. Model \eqref{eq:nv_dis_feature_dd1}--\eqref{eq:nv_dis_feature_dd2} is, in effect, fully data-driven.

Finally, to recast problem~\eqref{eq:nv_dis_feature_dd1}--\eqref{eq:nv_dis_feature_dd2} as a  linear program, we introduce the auxiliary variables $o_t$ and $u_t$ to equivalently reformulate the positive-part function as follows:
\begin{align}
    &\min_{\mathbf{q}, \mathbf{u}, \mathbf{o}} \quad \frac{1}{|\mathcal{T}|} \sum_{t\in\mathcal{T}} \psi_t^{-}u_t + \psi_t^{+}o_t\label{eq:nv_dis_feature_lp1}\\
    &\text{s. t.} \enskip u_t \geq \sum_{j=1}^p q^j x^j_t - E_t, \enskip \forall t\in\mathcal{T} \label{eq:nv_dis_feature_lp2}\\
    & \phantom{s. t.} \enskip o_t \geq E_t - \sum_{j=1}^p q^j x^j_t, \enskip \forall t\in\mathcal{T} \label{eq:nv_dis_feature_lp3}\\
    & \phantom{s. t.} \enskip 0 \leq \sum_{j=1}^p q^j x_t^j \leq \overline{E}, \enskip \forall t\in\mathcal{T} \label{eq:nv_dis_feature_lp4}\\
    & \phantom{s. t.} \enskip u_t, o_t \geq 0, \enskip \forall t\in\mathcal{T} \label{eq:nv_dis_feature_lp5}
\end{align}
Next, we explain how  we use the linear program~\eqref{eq:nv_dis_feature_lp1}--\eqref{eq:nv_dis_feature_lp5} to improve the tasks of renewable energy forecasting and trading.
\subsection{Renewable Energy Forecasting}\label{sec:forecasting}
Problem~\eqref{eq:nv_dis_feature_lp1}--\eqref{eq:nv_dis_feature_lp5} provides us with a simple, but effective procedure to enhance the quality of a given renewable energy forecast by exploiting auxiliary information. For this purpose, first we need to set $\psi_t^{-} = \psi_t^{+} = 1, \forall  t\in\mathcal{T}$, in~\eqref{eq:nv_dis_feature_lp1}--\eqref{eq:nv_dis_feature_lp5}. This results in the following linear programming problem:
\begin{align}
    &\min_{\mathbf{q}, \mathbf{u}, \mathbf{o}} \quad \frac{1}{|\mathcal{T}|} \sum_{t\in\mathcal{T}} u_t + o_t\label{eq:nv_dis_feature_for1}\\
    &\text{s. t.} \enskip u_t \geq \sum_{j=1}^p q^j x^j_t - E_t, \enskip \forall t\in\mathcal{T} \label{eq:nv_dis_feature_for2}\\
    & \phantom{s. t.} \enskip o_t \geq E_t - \sum_{j=1}^p q^j x^j_t, \enskip \forall t\in\mathcal{T}\label{eq:nv_dis_feature_for3}\\
    & \phantom{s. t.} \enskip 0 \leq \sum_{j=1}^p q^j x_t^j \leq \overline{E}, \enskip \forall t\in\mathcal{T}\label{eq:nv_dis_feature_for4}\\
    & \phantom{s. t.} \enskip u_t, o_t \geq 0, \enskip \forall t\in\mathcal{T}\label{eq:nv_dis_feature_for5}
\end{align}
where the coefficients $q^j$ of the linear decision rule are now optimized to learn the \emph{median} of the random renewable energy production. In other words, $\sum_{j=1}^p q^{*j} x^j_{t}$, with $q^{*j}$ being the optimal value of $q^{j}$ obtained from~\eqref{eq:nv_dis_feature_for1}--\eqref{eq:nv_dis_feature_for5}, is expected to be a good estimate of the median of the renewable energy production at time $t$.  This can be also inferred from the fact that model \eqref{eq:nv_dis_feature_for1}--\eqref{eq:nv_dis_feature_for5} minimizes the average absolute value of the renewable power deviations over the training sample, subject to the positiveness and capacity constraints \eqref{eq:nv_dis_feature_for4}.

Finally, we just have to include the renewable energy forecast we desire to improve as one of the regressors or features $x^j$ in the linear decision rule. The remaining features will then correspond to that extra information we want to take advantage of to enhance the quality of the renewable energy forecast. This extra information may be of a very different nature. For example, some of the features could correspond to categorical variables (hour of the day, day of the week ...) and others could be forecasts of potentially related stochastic variables. As a matter of fact, several features in vector $\mathbf{x}$ could represent forecasts on the renewable energy production of interest, but issued by different entities. The only condition for a piece of information to be treated as a feature is that it must be available at the time when the enhanced renewable energy forecast is to be generated.

In the particular application we present later on, we seek to improve the onshore wind power production forecast of the DK1 area of the pan-European electricity market that is issued every day by the Danish TSO. This benchmark is referred to as  BM throughout the rest of the paper. To this end, we use, as additional features, the forecasts of the wind power production in neighbouring regions that are produced by the respective TSOs in charge of those regions. We also introduce the constant feature $\mathbf{x}^1 = \mathbf{1}$ to correct for possible offsets.
The onshore DK1-wind power forecast issued by the Danish TSO (that is, model BM) is produced by the tool known as \emph{WindFor} \cite{WindFor}, a state-of-the-art software for forecasting wind power production at different scales that leverages numerical weather predictions (wind speed and direction), statistics and artificial intelligence.
\subsection{Renewable Energy Trading}\label{sec:trading}
In principle, model~\eqref{eq:nv_dis_feature_lp1}--\eqref{eq:nv_dis_feature_lp5} could be directly used for renewable energy trading without further ado. To this aim, we would just need to solve this problem for the enriched dataset $\left\{(E_t, \psi_t^{-}, \psi_t^{+},  \mathbf{x}_t), \forall t\in\mathcal{T} \right\}$ and thus, obtain the optimal coefficient vector $\mathbf{q}^*$ defining the linear decision rule (this is what we call \emph{model training}). Then, the bid $E^{D}_{t}$ to be submitted by the renewable energy producer to the day-ahead market for time period $t$ of the market horizon would be computed as
\begin{equation}
    E^{D}_{t} = \sum_{j=1}^p q^{*j} x^j_{t} \label{eq:bid_test}
\end{equation}
Unfortunately, we observe in practice that the direct application of model~\eqref{eq:nv_dis_feature_lp1}--\eqref{eq:nv_dis_feature_lp5} does not produce, in general, a bid more profitable than the expected-value bid (that is, the bid consisting in submitting the point forecast of renewable energy production to the day-ahead market). The reason for this has to do with the limited predictability of the marginal opportunity costs $\psi^{-}$ and $\psi^{+}$ (i.e., the absence of repeating patterns in the series of these costs). In effect, as shown in Fig.~2 of \cite{Skajaa2015}, the most sophisticated models for predicting $\psi^{-}$ and $\psi^{+}$ deliver forecasts that are completely uninformative or misleading for lead times beyond several hours into the future. However, the lead times required for partaking in the day-ahead market are usually longer than 12-14 hours. This empirical observation is, besides, supported by economic theory:  the balancing market price $\lambda^{B}$  represents a marginal cost for system imbalances in real time, which should be purely random. Consequently, the balancing market price should behave as a noise around the spot price $\lambda^{D}$. As a result, there is little in $\psi^{-}$ and $\psi^{+}$ that can be predicted for lead times longer than several hours. In this situation, the model flexibility introduced by the features in problem~\eqref{eq:nv_dis_feature_lp1}--\eqref{eq:nv_dis_feature_lp5} tends to produce overfitted linear decision rules, that is, rules that capture ``fictitious'' patterns of $\psi^{-}$ and $\psi^{+}$ in the historical/training dataset, but that do not repeat themselves beyond that set.

Against this background, in lieu of model~\eqref{eq:nv_dis_feature_lp1}--\eqref{eq:nv_dis_feature_lp5}, we propose to solve the following optimization problem:
\begin{align}
    &\min_{a, \mathbf{u}, \mathbf{o}} \quad \frac{1}{|\mathcal{T}|} \sum_{t\in\mathcal{T}} \psi_t^{-}u_t + \psi_t^{+}o_t\label{eq:nv_dis_trad1}\\
    &\text{s. t.} \enskip u_t \geq a \hat{w}_t - E_t, \enskip \forall t \in \mathcal{T} \label{eq:nv_dis_trad2}\\
    & \phantom{s. t.} \enskip o_t \geq E_t - a \hat{w}_t, \enskip \forall t \in \mathcal{T} \label{eq:nv_dis_trad3}\\
    & \phantom{s. t.} \enskip u_t, o_t \geq 0, \enskip \forall t \in \mathcal{T} \label{eq:nv_dis_trad4}
\end{align}
where the single feature of this model, namely, $\hat{w}$, represents the improved renewable energy forecast obtained from model~\eqref{eq:nv_dis_feature_for1}--\eqref{eq:nv_dis_feature_for5}. What we suggest for renewable energy trading is, therefore, a two-step procedure in which we first improve the renewable energy forecast by way of~\eqref{eq:nv_dis_feature_for1}--\eqref{eq:nv_dis_feature_for5} and then we correct such a forecast for trading by means of the substantially less flexible model~\eqref{eq:nv_dis_trad1}--\eqref{eq:nv_dis_trad4}.

As reported in \cite{HUBER2019}, in newsvendor problems (similar to the renewable energy trading problem we address here), the bulk of the economic gains we attain from data-driven procedures are linked to the improvement of the estimate of $E$ that we get. Following this rationale, we first use~\eqref{eq:nv_dis_feature_for1}--\eqref{eq:nv_dis_feature_for5}  to enhance such a estimate as much as possible, and then employ~\eqref{eq:nv_dis_trad1}--\eqref{eq:nv_dis_trad4} to account for mid-term patterns of $\psi^{-}$ and $\psi^{+}$ (the little that we can explain about these costs) in the market bid. Therefore, we compute the bid to be submitted to the day-ahead market for time period $t$ as
\begin{equation}
    E^{D}_{t} = a^{*} \hat{w}_{t} \label{eq:bid_test_2step}
\end{equation}
with $a^{*}$ being the optimal decision-rule coefficient delivered by~\eqref{eq:nv_dis_trad1}--\eqref{eq:nv_dis_trad4}.

In the following section, we elaborate on the application of this two-step procedure on a real experiment.

\section{Experiment Design and Model Training}
\label{sec:design&training}
Next, we describe in detail the experiment conducted to assess the performance of the data-driven models introduced in Sections~\ref{sec:forecasting} and~\ref{sec:trading} for renewable energy forecasting and trading, respectively. As previously mentioned, we focus on the \emph{onshore wind power produced in the DK1} area of the pan-European electricity market.

This section is divided in three parts. In the first one, we present the data gathered and the different models trained and tested. In the second and third parts, we introduce the metrics  used to quantify the performance of those models and elaborate on how we train them, in that order.

\subsection{Data and Features}

All the data employed in this research span from 01/08/2015 to 04/22/2019 and are fully published and freely available for download from the website of the Danish TSO \cite{energinet} and the ENTSO-e Transparency Platform \cite{entsoe}, which facilitates the reproducibility of our results. These data pertain to various features that either relate to the hour of the day and day of the week, or to day-ahead predictions about a number of potentially relevant variables, specifically, the total load, scheduled generation and solar power production in DK1, and wind power productions (onshore, offshore or both) in market zones adjacent to DK1, namely, zone 2 of Denmark (DK2), zone 2 of Norway (NO2), zones 3 and 4 of Sweden (SE3 and SE4, respectively), and the bidding zone of Germany, Austria and Luxembourg (DE-AT-LU). According to the Manual of Procedures (MoP) of the ENTSO-e Transparency Platform \cite{EntsoeMoP3.1}, these predictions should be made available in the platform by the different TSOs no later than 18:00 h of day ${\rm D-1}$ and span the 24 hours of the following day~D. However, some TSOs are temporarily failing to faithfully comply with the ENTSO-e's MoP. This is the case, for example, of the Danish TSO, which is generally uploading the day-ahead forecasts pertaining to DK1 and DK2  several hours late (in the early morning of day~D). Besides, the day-ahead forecasts are not accompanied with their issuance time stamp, which makes it impossible to determine  the exact time point in day D at which those forecasts were generated.  This implies that we cannot guarantee that all the forecasts we use as features in our models below are time-consistent, that is, that have been issued at the same time. We explain later on, right after the presentation of our models, how we deal with this time-consistency issue in order to guarantee a rigorous analysis and evaluation of the proposed approach.

 On a different front, the categorical features named as ``hour of the day'' and ``day of the week'' each comprise a group of 0/1 time series, specifically, 24 time series for the case of ``hour of the day'' and seven for ``day of the week.'' These series, besides, take on a value of one for all the time periods that correspond to the label of the feature, and zero otherwise. For example, for every hour of ``Monday'', only one of the seven series of the feature ``day of the week'' takes on the value one, whereas the value of the other six is set to zero.

%These data pertain to various features that either relate to categorical information, specifically, hour of the day and day of the week, or to predictions about a number of potentially relevant variables. These predictions, in turn, are issued by different TSOs and are available at a certain time point in day ${\rm D-1}$ for the 24 hours of the following day~D. The variables to which these day-ahead predictions refer correspond to the wind power productions (onshore, offshore or both) in market zones adjacent to DK1, namely, zone 2 of Denmark (DK2), zone 2 of Norway (NO2), zones 3 and 4 of Sweden (SE3 and SE4, respectively), and the bidding zone of Germany, Austria and Luxembourg (DE-AT-LU). The data also include day-ahead predictions for the total load, scheduled generation and solar power production in DK1.

We build and train  seven models of the type of~\eqref{eq:nv_dis_feature_for1}--\eqref{eq:nv_dis_feature_for5}. The first five of these models differ from one another by the number of features they exploit. More precisely,
\begin{description}%[left= 0pt]
\item \emph{Forecasting Model 1} (FM1), which only includes the day-ahead predictions of the on- and offshore wind power production in DK1.
\item \emph{Forecasting Model 2} (FM2), which results from adding the categorical variables ``hour of the day'' and ``day of the week'', and the day-ahead forecasts of solar power production, scheduled generation and total load in DK1 to model FM1.
\item \emph{Forecasting Model 3} (FM3), which is derived from model FM1 by adding the day-ahead forecasts of the on- and offshore wind power production in DK2.
\item \emph{Forecasting Model 4} (FM4), which results from model FM3 by adding the day-ahead forecasts of the onshore wind power production in NO2, DE-AT-LU, SE3 and SE4, and the day-ahead forecasts of offshore wind power production in DE-AT-LU.
\item \emph{Forecasting Model 5} (FM5), which includes all the previous features.
\item \emph{Utopian Model 1} (UM1), which is analogous to FM4, but uses the realized values of all the features that represent forecasts (that is, the actual outcomes of the associated stochastic processes), except, logically, for the DK1-onshore wind power forecast, which is what we seek to improve.
\item  \emph{Utopian Model 2 (UM2)}, which is also similar to model FM4, but uses the realized values of wind power production in SE3, SE4, NO2 and DE-AT-LU.

\end{description}

Models  UM1 and  UM2 are unrealizable in practice, as they assume perfect information on some of the features (instead of forecasts). However, their analysis is worthwhile  to assess the impact of the time-consistency issue we mentioned before, as we explain next. Models FM1-FM5 can be divided into two groups, namely:
\begin{enumerate}
    \item Models FM1, FM2 and FM3 only use information relative to DK1 and/or DK2, and therefore, we can ensure that they exploit time-consistent information as they are issued and uploaded to the ENTSO-e Transaparency Platform by the same entity, i.e., the Danish TSO.
    \item Models FM4 and FM5 also make use of information relative to the rest of bidding zones. Hence, we cannot guarantee that these models employ time-consistent information. However, the performance comparison between models FM3, FM4 and UM2 allows us to measure the impact of this possible time inconsistency. In fact, this impact is concluded to be negligible in Section \ref{sec:results_forecasting} (around 0.25-0.30 percentage points in terms of prediction performance).
\end{enumerate}
Finally, model UM1 provides us with an upper bound on how much the DK1-onshore wind power forecast issued by the Danish TSO could be improved with our methodology by enhancing the info on the features. More precisely, it allows us to quantify how much we would gain in prediction performance if we could use the actual realized values of the features that our models exploit instead of their forecasts.

%Models  UM1 and  UM2 are unrealizable in practice, as they assume perfect information on some of the features (instead of forecasts). However, their analysis is worthwhile. Indeed, model  UM1 provides us with an upper bound on how much the DK1-onshore wind power forecast issued by the Danish TSO could be improved with our methodology by enhancing the info on the features, while model UM2 allows us to further assess the contribution of the features related to NO2, DE-AT-LU, SE3 and SE4 to that improvement.

% We build and train four models of the type of~\eqref{eq:nv_dis_feature_for1}--\eqref{eq:nv_dis_feature_for5}. These models differ from one another by the number of features they exploit. More precisely,
% \begin{description}%[left= 0pt]
% \item \emph{Model 1} (M1), which only includes the day-ahead predictions of the on- and offshore wind power production in DK1.
% \item \emph{Model 2} (M2), which results from adding the categorical variables ``hour of the day'' and ``day of the week'', and the day-ahead forecasts of solar power production, scheduled generation and total load in DK1 to model M1.
% \item \textcolor{blue}{M4}
% \item \emph{Model 3} (M3), which also results from model M1, but in this case adding the day-ahead forecasts of the onshore wind power production in DK2, NO2, DE-AT-LU, SE3 and SE4, and the day-ahead forecasts of offshore wind power production in DK2 and DE-AT-LU.
% \item \emph{Model 4} (M4), which includes all the previous features.
% \item \textcolor{blue}{M4}
% \end{description}

As previously mentioned, the benchmark model BM we use for comparison and evaluation is the raw onshore DK1-wind power forecast issued by the Danish TSO and produced by the tool \emph{WindFor} [29]. Note that this forecast (that is, the output of \emph{WindFor}) is a feature (that is, an input) common to all the models listed above. This way, the ultimate goal of these models is to enhance the Danish TSO's forecast by exploiting the info carried by the other features considered. In selecting those other features that may be potentially most relevant to enhancing the onshore DK1-wind power forecast issued by the Danish TSO, we have limited ourselves to information that: i) pertains to DK1 and/or neighbouring bidding zones and ii) is published either in the ENTSO-e Transparency Platform or on Energinet.dk's website.

% For trading the onshore DK1-wind power production in the pan-European day-ahead market, we construct and train a fifth model (M5) of the type of~\eqref{eq:nv_dis_trad1}--\eqref{eq:nv_dis_trad4} that receives as input the wind power forecast $\hat{w}_t$ from model M3, which, as discussed later, is the model that exhibits the best overall prediction performance over the test set.

For trading the onshore DK1-wind power production in the pan-European day-ahead market, we construct an eighth model TM of the type of~\eqref{eq:nv_dis_trad1}--\eqref{eq:nv_dis_trad4} that receives as input the wind power forecast $\hat{w}_t$ from model  FM3, which, as discussed later, is   the simplest among the proposed models for wind power forecasting that exhibit the best overall prediction performance over the test set. Furthermore, in Section \ref{sec:results_trading}, we compare the market performance of the bid produced by model TM with that of the trading strategies consisting in directly bidding the point forecast given by models BM and FM3 into the day-ahead market. Even though the day-ahead market of Nordpool closes at 12:00h of D-1, that is, six hours before the deadline established by ENTSO-e to publish the information that our model TM uses, this is a technical issue (imposed by the source of data we are working with) that does not invalidate our analysis in Section \ref{sec:results_trading}, because models BM, FM3, and TM are here compared under the same conditions.

The wind power forecasts for the market zone DE-AT-LU are available on a 15-min time resolution, while the rest are given in hourly resolution. Consequently, we compute the hourly average values of the DE-AT-LU data series. Besides, some of the series have missing values, although the proportion of gaps in the data series relative to their length is negligible. We fill these gaps with a linear interpolation of the values in their extremes. Last but not least, in models   FM1--FM5, UM1 and UM2, every non-categorical feature is dynamically scaled by the maximum value of the feature that is observed in the training dataset. The target variable, that is, the onshore wind power production in DK1 is also scaled by the  most updated value of the wind power capacity installed in that zone that is available in \cite{entsoe}, which is  3669 MW. For convenience, all the data series are labelled using Coordinated Universal Time (UTC), which is also the time reference we use for our experiments.
% 2966 vs 3669
\subsection{Performance Metrics}
To evaluate the performance of the various forecasting models stemming from~\eqref{eq:nv_dis_feature_for1}--\eqref{eq:nv_dis_feature_for5}, we use the \emph{Mean Absolute Error} (MAE) and the \emph{Root Mean Square Error} (RMSE), i.e.,
\begin{align} \label{eq:mae_metric}
    {\rm MAE}:= \frac{1}{|\widetilde{\mathcal{T}}|} \sum_{t \in \widetilde{\mathcal{T}}} |E_{t} - E^D_{t}|
\end{align}
\begin{align} \label{eq:mnrmse_metric}
    {\rm RMSE}:= \frac{1}{|\widetilde{\mathcal{T}}|} \sqrt{\sum_{t \in \widetilde{\mathcal{T}}} (E_{t} - E^D_{t})^2}
\end{align}
where $\widetilde{\mathcal{T}}$ is the test set.

Recall that, when forecasting, the purpose of model~\eqref{eq:nv_dis_feature_for1}--\eqref{eq:nv_dis_feature_for5} is to improve an existing renewable energy prediction. In our case, this prediction is the day-ahead forecast of the onshore wind power production in DK1 that is issued by the Danish TSO every day. For this reason, we are especially interested in the percentage improvement with respect to that forecast in terms of MAE and RMSE.

On the other hand, to assess the performance of the trading model that results from~\eqref{eq:nv_dis_trad1}--\eqref{eq:nv_dis_trad4}, we compute the \emph{average opportunity loss} (AOL) linked to the onshore wind power production in DK1 over the test set, that is:
\begin{align} \label{eq:nv_metric}
    {\rm AOL}:= \frac{1}{|\widetilde{\mathcal{T}}|} \sum_{t \in \widetilde{\mathcal{T}}} \psi^{-}_{t} (E_{t} - E^D_{t})^+{+} \psi^{+}_{t} (E^D_{t} - E_{t})^+
\end{align}
The AOL gives us an idea of the monetary value lost by the onshore wind power production in DK1 due to its limited predictability. Therefore, rather than in the value of AOL \emph{per se}, we are interested in the decrease in AOL that we attain by means of model~\eqref{eq:nv_dis_trad1}--\eqref{eq:nv_dis_trad4} relative to the AOL delivered by submitting the Danish TSO's forecast to the day-ahead market.

Finally, note that if $\psi^{-}_{t}$ and $\psi^{+}_{t}$ are set to one for all $t$ the AOL metric becomes equivalent to computing the MAE.
\subsection{Model Training}
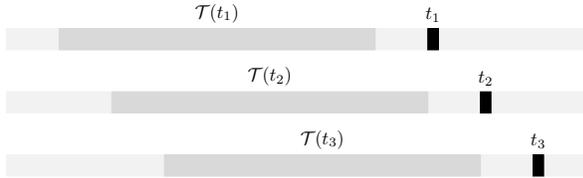
\begin{figure}
    \centering
    \begin{tikzpicture}[scale=0.7,every node/.style={scale=0.7}]
        \draw[draw=gray!10,fill=gray!10] (0,0) rectangle ++(11,0.4);
        \draw[draw=gray!30,fill=gray!30] (1,0) rectangle ++(6,0.4) node[above,pos=0.5,yshift=2mm] {$\mathcal{T}(t_1)$};
        \draw[draw=black,fill=black] (8,0) rectangle ++(0.2,0.4) node[above,pos=0.5,yshift=2mm] {$t_1$};
        \draw[draw=gray!10,fill=gray!10] (0,-1.2) rectangle ++(11,0.4);
        \draw[draw=gray!30,fill=gray!30] (2,-1.2) rectangle ++(6,0.4) node[above,pos=0.5,yshift=2mm] {$\mathcal{T}(t_2)$};
        \draw[draw=black,fill=black] (9,-1.2) rectangle ++(0.2,0.4) node[above,pos=0.5,yshift=2mm] {$t_2$};
        \draw[draw=gray!10,fill=gray!10] (0,-2.4) rectangle ++(11,0.4);
        \draw[draw=gray!30,fill=gray!30] (3,-2.4) rectangle ++(6,0.4) node[above,pos=0.5,yshift=2mm] {$\mathcal{T}(t_3)$};
        \draw[draw=black,fill=black] (10,-2.4) rectangle ++(0.2,0.4) node[above,pos=0.5,yshift=2mm] {$t_3$};
    \end{tikzpicture}
    \caption{Illustration of the rolling-window approach.}
    \label{fig:rolling}
\end{figure}

Except for the categorical information ``hour of the day'' and ``day of the week'', all the features we exploit in models~\eqref{eq:nv_dis_feature_for1}--\eqref{eq:nv_dis_feature_for5} are forecasts of a variety of potentially informative variables for time $t$. All these forecasts pertain to the 24 hours of the following day. In actual practice, models~\eqref{eq:nv_dis_feature_for1}--\eqref{eq:nv_dis_feature_for5} and~\eqref{eq:nv_dis_trad1}--\eqref{eq:nv_dis_trad4} are trained using a rolling-window approach and therefore, the training set depends on each time period $t$ of the test set $\widetilde{\mathcal{T}}$. The rolling training set is denoted here as $\mathcal{T}(t)$ and is illustrated in Fig. \ref{fig:rolling}. Notice that the length of the training set is kept constant as time progresses. Furthermore, there exists a gap between the time period $t$ and its corresponding training set $\mathcal{T}(t)$. The reason for this gap is that the values of $E$, $\psi^{-}$ and $\psi^{+}$ for the time interval that goes from the moment the forecasts are made available and time $t$ are still not known and consequently, such time periods cannot be used for the training of the  models~\eqref{eq:nv_dis_feature_for1}--\eqref{eq:nv_dis_feature_for5} or~\eqref{eq:nv_dis_trad1}--\eqref{eq:nv_dis_trad4}.
%Therefore, from the training set, we must remove the hours in the time interval that goes after the moment our forecast is supposed to be issue and time $t$.
%
%The wind power producer must place a bid into the day-ahead market at a certain time period $t'$ of day ${\rm D-1}$ (for example, at noon in DK1) for the 24 hours of the following day D. Then, at that instant $t'$, the training set includes all the past data on the features from and including, hour 23:00 (UTC) of day ${\rm D-2}$ back to hour 00:00 (UTC) of day ${\rm D-1}-\frac{T}{24}$. For simplicity, we omit all the  hours of day ${\rm D-1}$ prior to time $t'$ from the training set. Note, however, that the loss of information implied by this omission is negligible, because in models~\eqref{eq:nv_dis_feature_for1}--\eqref{eq:nv_dis_feature_for5} and~\eqref{eq:nv_dis_trad1}--\eqref{eq:nv_dis_trad4} the time dynamics of the response variables (for instance, their autorregresive nature) is captured through the dynamics of the underlying explanatory features.
%
%Recall that, except for the categorical information ``hour of the day'' and ``day of the week'', all the features we exploit in models~\eqref{eq:nv_dis_feature_for1}--\eqref{eq:nv_dis_feature_for5} are predictions for a variety of potentially informative variables over the 24 hours of day ${\rm D}$. We assume, therefore, that at the moment of bidding into the day-ahead market, i.e., at time $t'$, these predictions are available, which is typically the case.

This rolling-window approach allows us to dynamically re-estimate the decision-rule parameters $\mathbf q$ and $a$ in~\eqref{eq:nv_dis_feature_for1}--\eqref{eq:nv_dis_feature_for5} and~\eqref{eq:nv_dis_trad1}--\eqref{eq:nv_dis_trad4}, respectively, as the information on the considered features is updated. Every time these parameters are re-estimated, equations \eqref{eq:bid_test} and \eqref{eq:bid_test_2step} are used  to issue improved forecasts and bids for time period $t$.

%The test set is made up of ``days {\rm D}'', where both the 24-h wind power forecast delivered by model~\eqref{eq:nv_dis_feature_for1}--\eqref{eq:nv_dis_feature_for5} and the 24 wind power bids given by model~\eqref{eq:nv_dis_trad1}--\eqref{eq:nv_dis_trad4} are evaluated in terms of MAE and AOL, respectively. In other words, every time the decision-rule parameters $\mathbf q$ and $a$ are re-estimated using the training set linked to day ${\rm D-1}$, these decision rules are employed to issue forecasts and bids for day $D$. Fig. \ref{fig:rolling_window} clarifies the rolling-window approach we follow.
%
Critical to the training of models FM1--FM5 and TM is determining the length $|\mathcal{T}(t)|$ of the training set. This length defines when the data linked to certain days in the past have become too old to be considered in the training process. We devote the first year of data to tune this length for models   FM1--FM5. In this time interval, the piece of data spanning from 08/07/2015 to 02/02/2016 (180 days) is used as the validation subset. We then compute the MAE on this subset for each of the models   FM1--FM5 and for different lengths of the training subset, which we vary from one to seven months.  We remark that the length of the training set is the only hyper-parameter that needs to be tuned for our models, which represents an advantage in terms of ease of use and implementation.

%\begin{table}[h!]
%\centering
%\begin{tabular}{c  c c c c }
%\hline
%month & M1         & M2         & M3         & M4         \\
%\hline
%1     & 0.046664 & 0.048919 & 0.050413 & 0.053929 \\
%2     & 0.046333 & 0.047036 & 0.047453 & 0.048611 \\
%3     & 0.046076 & 0.046804 & 0.046161 & 0.047087 \\
%4     & 0.046220 & 0.046729 & 0.046096 & 0.046744 \\
%5     & 0.046248 & 0.046962 & 0.045933 & 0.046486 \\
%6     & 0.046135 & 0.046623 & 0.045934 & 0.046127 \\
%7     & 0.046245 & 0.046565 & 0.045945 & 0.046295 \\
%\hline
%\end{tabular}
%\caption{Evolution of MAE as per unit of the DK1 wind power capacity for different lengths of the training set. {\color{blue} DAR OK}}
%\label{table:mae_training_lenght_unit}
%\end{table}

\begin{table}[h!]
\centering
\begin{tabular}{c c c c c c}
\hline
$|\mathcal{T}(t)|$ &  FM1 &  FM2 &  FM3  &  FM4 &  FM5     \\
\hline
1     & 11.67 & 7.40  &  10.37   & 4.57  & -2.08 \\
2     & 12.30 & 10.97 &  11.95   & 10.18 & 7.98  \\
3     & 12.78 & 11.40 &  12.55   & 12.62 & 10.87 \\
4     & 12.51 & 11.55 &  12.38   & 12.75 & 11.52 \\
5     & 12.46 & 11.10 &  12.48   & 13.05 & 12.01 \\
6     & \textbf{12.67} & \textbf{11.75} &   \textbf{12.83} &  \textbf{13.05} & \textbf{12.69} \\
7     & 12.46 & 11.86 &  12.49  & 13.03 & 12.37 \\
\hline
\end{tabular}
\caption{MAE reduction in percentage with respect to the benchmark for different lengths of the training set in months.}
\label{table:mae_training_lenght_perc}
\end{table}

Table~\ref{table:mae_training_lenght_perc} summarizes the results of this analysis, where the MAE linked to each model and length is expressed in percentage reduction with respect to the MAE associated with the benchmark, namely, the onshore DK1-wind power forecast issued by the Danish TSO.  From this table, we notice that the improvement in the performance of models FM1--FM5 that is initially observed as we increase the length of the training set not only ends up saturating, but even reverses as we extend the training set beyond several months (e.g., six months in the case of FM5). This is due to the fact that, at some point in time, the information contained in the oldest data becomes obsolete and thus, potentially misleading. In light of these results, we set the length of the training set for forecasting to six months.

We proceed in a similar fashion to establish the length of the dataset we use to train the trading model TM. In this case, we change the validation subset to 06/03/2016--11/29/2016 (180 days). This change is required because model TM is fed with the improved wind power forecast yielded by FM3  (the one exhibiting the best trade-off between simplicity, data reliability and forecasting performance). Consequently, training model TM involves generating a sufficient number of predictions from model FM3 first, which, in turn, is to be trained over a dataset spanning six months. Hence, we need to reserve a big chunk of data to study the impact of the length of the training set on the performance of model TM. Table~\ref{table:aol_training_lenght} shows the results of this study for a length of the training set varying from one to ten months. The numbers in the table correspond to the AOL reduction of model TM in percentage with respect to the AOL given by the benchmark, that is, the trading strategy consisting in submitting the wind power prediction issued by the Danish TSO to the pan-European electricity market. In view of these results and for ease of implementation, we also set the length of the training set for trading to six months.

\begin{table}[h!]
\setlength{\tabcolsep}{3pt}
\centering
\begin{tabular}{c c c c c c c c c c c}
% \hline
% month & M5   & M5(\%) \\
% \hline
% 1     & 0.087338   & 2.51\%  \\
% 2     & 0.082906   & 7.45\%  \\
% 3     & 0.082980   & 7.37\%  \\
% 4     & 0.083232   & 7.09\%  \\
% 5     & 0.083533   & 6.75\%  \\
% 6     & \textbf{0.082225}   & \textbf{8.21}\%  \\
% 7     & 0.083533   & 6.75\%  \\
% 8     & 0.083931   & 6.31\%  \\
% 9     & 0.084236   & 5.97\%  \\
% 10    & 0.084533   & 5.64\%  \\
% \hline
\hline
$|\mathcal{T}(t)|$ & 1 & 2 & 3 & 4 & 5 & 6 & 7 & 8 & 9 & 10 \\
%M5 & 2.51 & 7.45 & 7.37 & 7.09 & 6.75 & \textbf{8.21} & 6.75 & 6.31 & 5.97 & 5.64 \\
  TM &   3.04 &   7.62 &   8.20 &   8.33 & \textbf{9.29} & \textbf{9.14} &   8.36 &   8.25 &   8.05 &   7.67 \\
\hline
\end{tabular}
\caption{AOL reduction in percentage of model TM with respect to the benchmark for different lengths of the training set in months.}
\label{table:aol_training_lenght}
\end{table}

Next, we discuss the results obtained from the simulation conducted on all the remaining days in the full dataset that have not been used to determine the length of the training set.

\section{Results}
\label{sec:results}

We divide this section in two parts. In the first one, we present and discuss the improvements in wind power forecasting brought about by the linear decision rule that results from~\eqref{eq:nv_dis_feature_for1}--\eqref{eq:nv_dis_feature_for5}. Subsequently, we elaborate on the improvements in wind power trading that we attain by means of model~\eqref{eq:nv_dis_trad1}--\eqref{eq:nv_dis_trad4}.
\subsection{Improvements in Wind Power Forecasting}
\label{sec:results_forecasting}
The first and last days in the test set are 02/04/2016 and 04/22/2019. That is, the test set in the simulation comprises 1174 days in total.
Table \ref{table:mae_rmse_improvement} provides the MAE and the RMSE reductions (in percentage) with respect to the performance metrics of the benchmark, namely, the raw forecast issued by the Danish TSO.

%{\color{blue} The difference between the reference value (prediction of the Danish TSO) and the value of the performance metric for each model and each metric is obtained. Then, the percentage improvements showed in table \ref{table:mae_rmse_improvement} are computed as the quotient between this difference and the mentioned reference value.}
%
\begin{table}[h!]
\centering
 \begin{tabular}{c c c c c c c c}
 \hline
  &  FM1 &  FM2 &  FM3 &  FM4 &  FM5 &  UM1 &  UM2 \\ [0.5ex]
 \hline
%  \multirow{2}{3em}{MAE} & 0.037024 & 0.034423 & 0.034421 & 0.033859 & 0.033865 \\
MAE & 7.03 & 7.03 &  8.53 & 8.55 & 8.53 &  10.18 &  8.80 \\
%  \hline
%  \multirow{2}{3em}{RMSE} & 0.046578 &	0.043840 & 0.043749 & 0.043231 & 0.043176 \\
RMSE & 6.04 & 6.22 &  7.16 & 7.33 & 7.46 &  9.14 &  7.46 \\
 \hline
 \end{tabular}
\caption{MAE and RMSE reduction (in percentage) with respect to the benchmark.}
\label{table:mae_rmse_improvement}
\end{table}
We observe that model~\eqref{eq:nv_dis_feature_for1}--\eqref{eq:nv_dis_feature_for5}, which, in essence, is an easily implementable and computationally inexpensive linear program, is able to substantially enhance the wind power forecasts made by Energinet.dk. Actually, most of the reduction can be achieved by linearly combining Energinet.dk's predictions for the onshore and offshore DK1-wind power productions (model  FM1). From these results, we infer that historical information of wind power forecasts pertaining to neighboring bidding zones is not currently being exploited by the Danish TSO. In contrast, the performance comparison of models FM1 and FM2, on the one hand, and of models FM4 and FM5, on the other, seems to signal the fact that the ``hour of the day'' and ``day of the week'', and the day-ahead forecasts of solar power production, scheduled generation and total load in DK1 do not have predictive power on the targeted variable. Besides, the comparison between FM1 and UM1 reveals that a significant improvement in the forecasting of the DK1-onshore wind power production can be attained by enhancing the quality of the forecasts in the areas adjacent to DK1. However, the comparison between FM1, FM3, FM4 and UM2 further indicates that the bulk of this potential improvement is to be attributed to the DK1 and DK2 wind-related features. The models of the type of~\eqref{eq:nv_dis_feature_for1}--\eqref{eq:nv_dis_feature_for5} that exhibit the best forecasting performance are FM3, FM4 and FM5. Since the first one is significantly simpler than the rest and offers the best guarantees in terms of data reliability\footnote{Model FM3 uses data produced and uploaded at the same time to the ENTSO-e Transparency Platform by the same entity, namely, the Danish TSO. One can expect, therefore, that the forecasts contained in these data have been built with the same past information available.}, we use FM3 to feed TM with the required wind power forecast. Interestingly, even though models FM2 and FM5 exploit a larger number of features than FM1 and FM4, respectively, their forecasting performance is not (or barely) improved.

\newcommand{\intervalFigSampleForecast}{01/01/16 to 01/08/16}

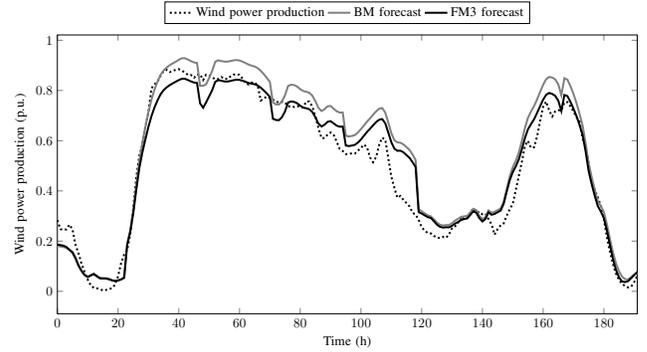
\begin{figure}
\centering
\begin{tikzpicture}[scale=0.5]
	\begin{axis}[	
    width=17cm, %\textwidth,
    height=9cm,
    xmin = 0,
    xmax = 191,
	legend style={at={(0.5,1.04)},anchor=south,legend cell align=left,legend columns=3},
	clip marker paths=true,	
	xlabel = Time (h),
	ylabel = Wind power production (p.u.)]
	\addplot[line width=1.5pt,draw=black,dotted] table [x=index, y=real, col sep=comma] {data_prediction_results_fm3.csv}; \addlegendentry{Wind power production}
	\addplot[line width=1.5pt,draw=gray] table [x=index, y=energinet, col sep=comma] {data_prediction_results_fm3.csv}; \addlegendentry{BM forecast}
    \addplot[line width=1.5pt,draw=black] table [x=index, y=FM3, col sep=comma] {data_prediction_results_fm3.csv}; \addlegendentry{FM3 forecast}
	\end{axis}	
\end{tikzpicture} \\
\caption{Illustration of the forecasts issued by the Danish TSO ( BM) and model  FM3 for the interval \intervalFigSampleForecast.}   \label{fig:sample_forecast}
\end{figure}

For the sake of illustration, Fig. \ref{fig:sample_forecast} plots the actual realization of the wind power production in the time interval \intervalFigSampleForecast, together with the forecasts issued by Energinet.dk (BM) and the proposed model FM3. It can be observed that from hour 80 on, the forecast yielded by FM3 is always closer to the actual wind power production than the forecast used by the Danish TSO. On average, model FM3 produces forecasts that, over the simulation period, deviate  100.44 MW with respect to the true wind power values, whereas Energinet.dk's average deviation for this period amounts to 109.82 MW.

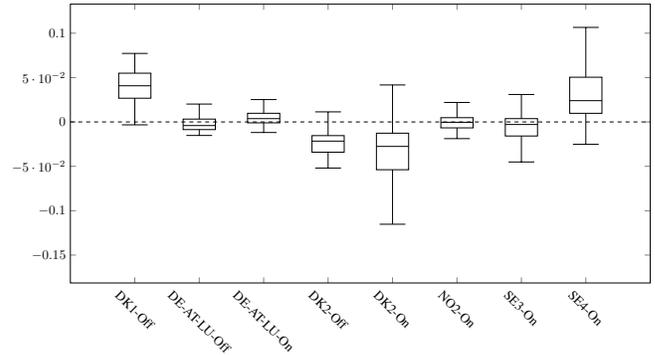
\begin{figure}
\centering
\begin{tikzpicture}[scale=0.5]
\begin{axis}[
    width=17cm, %\textwidth,
    height=9cm,
    xmin = 0,
    xmax = 9,
	legend style={at={(0.5,1.04)},
	anchor=south,
	legend cell align=left,
	legend columns=3},	
% 	clip marker paths=true,	
    boxplot/draw direction=y,
    xtick={1,2,3,4,5,6,7,8,9},
    xticklabels={ DK1-Off, DE-AT-LU-Off, DE-AT-LU-On, DK2-Off, DK2-On, NO2-On, SE3-On, SE4-On},
    xticklabel style={rotate=-45},
    %ylabel = Coefficient values
    ]	
    %\addplot [line width=0.5pt, draw=black, boxplot={box extend=0.5}] table [y index=2, col sep=comma]{data_tabla_coeficientes_caso3.csv};
    \addplot[mark=none, dashed] coordinates {(0,0) (10,0)};
    \addplot [line width=0.5pt, draw=black, boxplot={box extend=0.5}] table [y index=4, col sep=comma]{data_tabla_coeficientes_caso3.csv};
    \addplot [line width=0.5pt, draw=black, boxplot={box extend=0.5}] table [y index=5, col sep=comma]{data_tabla_coeficientes_caso3.csv};
    \addplot [line width=0.5pt, draw=black, boxplot={box extend=0.5}] table [y index=6, col sep=comma]{data_tabla_coeficientes_caso3.csv};
    \addplot [line width=0.5pt, draw=black, boxplot={box extend=0.5}] table [y index=7, col sep=comma]{data_tabla_coeficientes_caso3.csv};
    \addplot [line width=0.5pt, draw=black, boxplot={box extend=0.5}] table [y index=8, col sep=comma]{data_tabla_coeficientes_caso3.csv};
    \addplot [line width=0.5pt, draw=black, boxplot={box extend=0.5}] table [y index=9, col sep=comma]{data_tabla_coeficientes_caso3.csv};
    \addplot [line width=0.5pt, draw=black, boxplot={box extend=0.5}] table [y index=10, col sep=comma]{data_tabla_coeficientes_caso3.csv};
    \addplot [line width=0.5pt, draw=black, boxplot={box extend=0.5}] table [y index=11, col sep=comma]{data_tabla_coeficientes_caso3.csv};

\end{axis}
\end{tikzpicture}
\caption{Box plot of the coefficients obtained for FM4 in the simulation period 02/04/16 to 04/22/2019.}
\label{fig:coeficient_boxplot}
\end{figure}

The simplicity of model~\eqref{eq:nv_dis_feature_for1}--\eqref{eq:nv_dis_feature_for5} makes it more interpretable than other forecasting models based, for instance, on artificial neural networks. Not surprisingly, the coefficient corresponding to the onshore DK1-wind power forecast issued by the Danish TSO is the largest one for the FM and UM models. For example, its value in model FM4 ranges from 0.8335 to 1.0267 over the simulation period. The other coefficient values of model FM4 are depicted in a box plot in Fig. \ref{fig:coeficient_boxplot}. As observed, the forecasts for the offshore DK1-wind, the onshore and offshore DK2-wind, and the onshore SE4-wind are also significant.
%{\color{blue}0.8335 to 1.0267} para FM4 (M3)
%{\color{blue}0.8273 to 1.0148} para FM3

%Table \ref{table:coeficient_statistics} includes {\color{red} the minimum and maximum values, and the quartiles} of the frequency distribution of the coefficient vector the average value taken on by the coefficient vector $\mathbf{q}$ computed over the simulation period for model variant M3. Not surprisingly, the largest coefficient corresponds to the onshore DK1-wind power forecast issued by the Danish TSO, albeit the contributions of the forecasts for the offshore DK1-wind, the onshore and offshore DK2-wind, and the onshore SE4-wind are also significant. % the wind power productions in neighbouring areas also show significant contributions.

% \begin{figure}
% \centering
% \begin{tikzpicture}[scale=0.5]
% \begin{axis}[
%     width=17cm, %\textwidth,
%     height=9cm,
% 	legend style={at={(0.5,1.04)},
% 	anchor=south,
% 	legend cell align=left,
% 	legend columns=3},	
% % 	clip marker paths=true,	
%     boxplot/draw direction=y,
%     xtick={1,2,3,4},
%     xticklabels={DK1-Off, DK2-On, DK2-Off},
%     xticklabel style={rotate=-45},
%     %ylabel = Coefficient values
%     ]	
%     %\addplot [line width=0.5pt, draw=black, boxplot={box extend=0.5}] table [y index=2, col sep=comma]{data_tabla_coeficientes_caso3.csv};
%     \addplot [line width=0.5pt, draw=black, boxplot={box extend=0.5}] table [y index=4, col sep=comma]{data_tabla_coeficientes_fm3.csv};
%     \addplot [line width=0.5pt, draw=black, boxplot={box extend=0.5}] table [y index=5, col sep=comma]{data_tabla_coeficientes_fm3.csv};
%     \addplot [line width=0.5pt, draw=black, boxplot={box extend=0.5}] table [y index=6, col sep=comma]{data_tabla_coeficientes_fm3.csv};
% \end{axis}
% \end{tikzpicture}
% \caption{Box plot of the coefficients obtained for {\color{blue}FM3} in the simulation period 02/04/16 to 04/22/2019.}
% \label{fig:coeficient_boxplot}
% \end{figure}

\subsection{Improvements in Wind Power Trading}
\label{sec:results_trading}
The first and last days of the test set, in this case, are 11/30/2016 and 04/22/2019, in that order. This means that the test set in this simulation consists of 874 days. In this analysis we assume that the wind power point forecast issued by each model is directly bid into the day-ahead market and then we compute the average opportunity loss as in \eqref{eq:nv_metric}.

If the forecasts issued by FM3 are used as bids, the AOL is reduced by 1.30\% with respect to the benchmark, which consists in bidding the raw wind power point forecast issued by the Danish TSO into the day-ahead market. Although model FM3 is tailored to forecasting, the reduction of the prediction error that it achieves is accompanied with an AOL decrease too. %Although model M3 is tailored to forecasting, the reduction of the prediction error that it achieves inevitably involves an AOL reduction too.

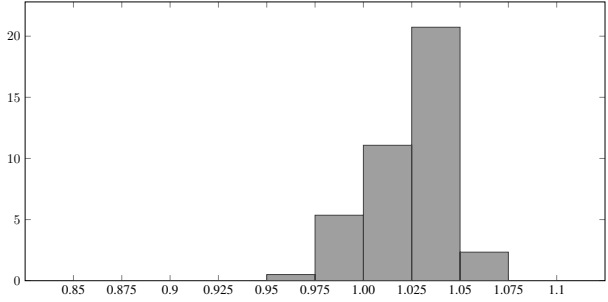
\begin{figure}
  \centering
    \begin{tikzpicture} [scale=0.5]
    \begin{axis}[
        width=17cm, %\textwidth,
        height=9cm,
        ymin=0,
        xtick={0.85, 0.875, 0.9, 0.925, 0.95, 0.975, 1.00, 1.025, 1.05, 1.075, 1.1},
        xticklabels={0.85, 0.875, 0.9, 0.925, 0.95, 0.975, 1.00, 1.025, 1.05, 1.075, 1.1},
        legend style={at={(0.5,1.04)},anchor=south,legend cell align=left,legend columns=3},	
	    clip marker paths=true,	
        ]
        \addplot [hist={data=y, data min=0.85, data max=1.1, density}, fill=gray!75,draw=gray!50!black] table
        [y index=2, col sep=comma]{data_coeficientes_seg_paso_fm3.csv};
        %\addlegendentry{coeficientes}
        %\addplot [line width=2.5pt, draw=black] coordinates{(1.025, 0) (1.025, 15)}; \addlegendentry{median}
    \end{axis}
    \end{tikzpicture}
    \caption{Histogram of the values taken on by the decision-rule parameter $a$ in model TM for the interval 11/30/16 to 04/22/19.}
  \label{fig:histograma}
\end{figure}

If the mid-term dynamics of the marginal opportunity costs are accounted for through model TM, the AOL reduction increases up to 2.13\%. In this regard, the histogram of the values taken on by the decision-rule parameter $a$ over the simulation period is plotted in Fig.~\ref{fig:histograma}. Interestingly, this parameter tends to take values above 1, so as to profit from the fact that, in the DK1 bidding zone, overproduction is, on average, more penalized than underproduction.

To further explain the AOL reduction achieved by TM, we define next the empirical \textit{critical fractile} estimated over the training set $\mathcal{T}$ as:
\begin{equation}
    R =  \frac{\frac{1}{|\mathcal{T}|}\sum_{t \in \mathcal{T}} \psi^+_t}{\frac{1}{|\mathcal{T}|}\sum_{t \in \mathcal{T}} \psi^-_t + \frac{1}{|\mathcal{T}|}\sum_{t \in \mathcal{T}} \psi^+_t}
\end{equation}
The ratio $R$ balances the marginal opportunity cost for overproduction and the marginal opportunity cost for either under- or overproduction, all of them averaged over  $\mathcal{T}$. A value of $R$ higher than 0.5 means that the opportunity cost for overproduction was more significant than that for underproduction throughout the training period. In such a case, the optimal market bid should be higher than the forecast production in order to hedge against overproduction. Conversely, if $R$ is lower than 0.5, the optimal market bid should be lower than the forecast production.

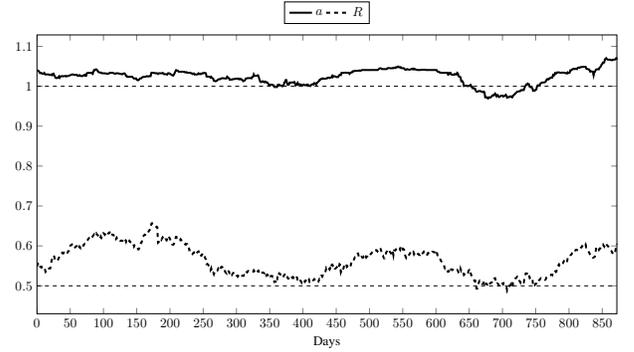
\begin{figure}
\centering
\begin{tikzpicture}[scale=0.5]
	\begin{axis}[	
    width=17cm, %\textwidth,
    height=9cm,
    xmin = 0,
    xmax = 872,
	legend style={at={(0.5,1.04)},anchor=south,legend cell align=left,legend columns=3},	
	clip marker paths=true,	
	xlabel = Days]
	\addplot[line width=1.5pt,draw=black] table [x=index, y=var_a, col sep=comma] {data_tendencia_fm3.csv}; \addlegendentry{$a$}
	\addplot[line width=1.5pt,draw=black,dashed] table [x=index, y=r_trn, col sep=comma] {data_tendencia_fm3.csv}; \addlegendentry{$R$}
	\addplot[mark=none, dashed] coordinates {(0,0.5) (870,0.5)};
	\addplot[mark=none, dashed] coordinates {(0,1) (870,1)};
	\end{axis}	
\end{tikzpicture} \\
\caption{Evolution of decision-rule parameter $a$ in TM and ratio $R$ for the interval 11/30/16 to 04/22/19.}
\label{fig:tendency}
\end{figure}

\begin{figure}
\centering
\begin{tikzpicture}[scale=0.5]
	\begin{axis}[	
    width=17cm, %\textwidth,
    height=9cm,
    xmin = 0,
    xmax = 873,
	legend style={at={(0.5,1.04)},anchor=south,legend cell align=left,legend columns=3},	
	clip marker paths=true,	
	xlabel = Days,
	ylabel = Euro]
	\addplot[line width=1.5pt,draw=black] table [x=index, y=accum3669x24, col sep=comma] {data_aol_accumulated.csv};
	\end{axis}	
\end{tikzpicture} \\
\caption{Accumulated opportunity-loss reduction of TM for the interval 11/30/16 to 04/22/19 for a installed capacity of 3669 MW.}
\label{fig:accumulated_aol}
\end{figure}
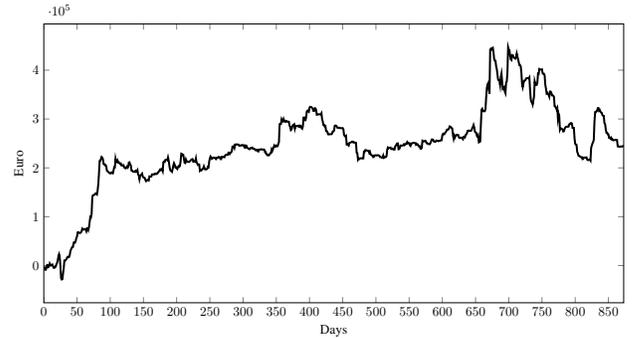

Fig. \ref{fig:tendency} depicts the time evolution of the decision-rule parameter $a$ in TM together with the ratio $R$ over the simulation period 11/30/16-04/22/19. As observed, the value of $a$ continuously adapts to the variations of $R$ as the training period $\mathcal{T}$ moves forward. This way, the bids provided by TM take into account the mid-term dynamics of $\psi^{-}$ and $\psi^{+}$ to properly hedge against under or overproduction.
Finally,  Fig. \ref{fig:accumulated_aol} illustrates the accrued reduction in opportunity loss achieved by model TM with respect to the benchmark over the simulation period. Note that the plot is studded with time instants when the accrued improvement suddenly decreases. This is because the series of balancing prices is scattered with highly unpredictable spikes. Indeed, the limited predictability of balancing prices is what makes the trading strategy consisting in minimizing expected deviations so hard to beat.  To finish this section, we note that a similar experiment could be conducted for a particular wind farm or a particular wind power producer. To illustrate the benefits of our approach, however, we have decided to work on the aggregate onshore DK1-wind power production for several reasons. First, there is a rich set of market data related to DK1 and sorrounding bidding zones publicly available in the ENTSO-e Transparency Platform, while analogous datasets for particular wind farms or producers are usually kept confidential. Second, DK1 uses a dual-price balancing settlement, and lastly, the onshore installed wind power capacity in DK1 amounts to 3669 MW. Today, portfolios of similar size can be easily found in countries such as Spain, United Kingdom or Germany \cite{WPDatabase}.
\section{Conclusions}
\label{sec:conclusions}
In this paper we have proposed an inexpensive, easy-to-implement, but effective method to enhance the tasks of renewable energy forecasting and trading. Our method is based on a data-driven newsvendor-type optimization model that leverages extra available information to produce an improved renewable energy forecast or a renewable energy bid that can be directly placed in the day-ahead electricity market.
%and auxiliary information to deliver both, an improved wind power forecast and a direct data-driven bidding.

The effectiveness of our approach is tested on a realistic case study where we aim, on the one hand, to improve the forecast issued by the Danish TSO for the onshore wind power production in the DK1 bidding zone of the pan-European electricity market, and, on the other, to formulate a competitive market bid for such a production. To this end, we build a rolling-window simulation setup that mimics the actual processes of forecasting and bidding and exploits the information available at the moment the forecast must be issued or the bid must be placed.
%in the forecasting and bidding process through a rolling window. %Fig \ref{fig:rolling_window} aims to clarify this setup.

The numerical results highlight the benefits achieved by our approach, which amounts to a 8.53\% of reduction in MAE and a 2.13\% of improvement of AOL with respect to the benchmarks for the simulation period considered. These figures point out the intrinsic value of exploiting additional information such as spatially correlated forecasts. In this line, we have observed that the use (as features) of both on- and offshore wind power forecasts in areas geographically close to the zone to which the target wind power production belongs are valuable. This seems to be especially true if those areas pertain to the same country or domain of the same TSO.

Future work could be focused on the development of robust counterparts of the proposed models with the aim of reducing the volatility of the improvements achieved. Variable selection methodologies could also be implemented to determine the best subset of regressors to feed in the models and to enhance model interpretability. Likewise, nonlinear mappings between the features and the response variable could be captured within our approach by performing nonlinear transformations on the features or by way of kernels. Neither of these actions would significantly affect the computational complexity of our models. More generally, adapting theory and techniques from the rich field of nonlinear regression to the data-driven framework we propose is indeed a relevant path to follow in future research.

%Future work could be focused on the development of robust counterparts of the proposed models with the aim of reducing the volatility of the improvements achieved. Variable selection methodologies could also be implemented as a previous step to determine the best subset of regressors to feed in the models. Likewise, nonlinear relationships between the features and the response variable could be easily captured within our approach by performing nonlinear transformations on the features or by way of kernels. Neither of these actions would significantly affect the computational complexity of our models, which would remain as linear programming problems.

Finally, while the data-driven model for renewable energy trading we develop is tailored to electricity markets with a dual-price settlement for imbalances, such as Nordpool-DK1 or MIBEL (Spain), it could also be adapted to any market where deviations with respect to a predefined forward schedule entail an opportunity cost.

% Table \ref{table:mae_training_lenght_unit} and \ref{table:aol_training_lenght} point out the crucial role of the length of the training set that, on behalf of this results, should be consider as any other hyper-parameter. Nevertheless, meanwhile in a forecast context (Table \ref{table:mae_training_lenght_unit}) the effect of the training length seems to saturate as it increases, this conclusion does not hold when bidding due to mid-term cost patterns.

% Table \ref{table:mae_rmse_improvement} shows an improvement of 8.55\% in MAE metric. Table \ref{table:mae_rmse_improvement} also enlighten the outstanding effect of considering spatial information, i.e., the forecast of neighbour areas, when it comes to improve the the prediction of a particular location.

% The per unit income surplus obtained with the strategy described in this manuscript could be vet in Fig. \ref{fig:accumulated_aol}. Although a clearly positive tendency is observed, attention should be paid to the inherent volatile short-term dynamic due to the random price behaviour.

\bibliographystyle{ieeetr}

% \bibliography{bib_bibliografia_review}

\end{document}